\documentclass[12pt]{article}
\pdfoutput=1
\usepackage{natbib}
\usepackage{microtype}
\usepackage{graphicx}
\usepackage{subfigure}
\usepackage{booktabs} 
\usepackage{fullpage}
\usepackage{multirow}

\usepackage[breaklinks,colorlinks,
            linkcolor=red,
            anchorcolor=blue,
            citecolor=blue
            ]{hyperref}

\usepackage{amsmath}
\usepackage{amssymb}
\usepackage{mathtools}
\usepackage{amsthm}
\usepackage[figuresright]{rotating}
\usepackage[capitalize,noabbrev]{cleveref}
\usepackage{url}
\usepackage{footmisc}
\theoremstyle{plain}

\theoremstyle{definition}

\theoremstyle{remark}

\usepackage[textsize=tiny]{todonotes}

\usepackage{relsize}
\usepackage[utf8]{inputenc} 
\usepackage[T1]{fontenc}    
\usepackage{url}            
\usepackage{booktabs}       
\usepackage{amsfonts}       
\usepackage{nicefrac}       
\usepackage{microtype}      
\usepackage{smile}
\usepackage{bbm}
\usepackage{xcolor}

\def \d {\text{d}}

\title{\LARGE SOLNP+: A Derivative-Free Solver for Constrained Nonlinear Optimization}
\author{Dongdong Ge\footnote{Shanghai University of Finance and Economics. Email: \texttt{ge.dongdong@mail.shufe.edu.cn}}\quad   Tianhao Liu\footnote{Shanghai University of Finance and Economics. Email: \texttt{liu.tianhao@163.sufe.edu.cn}}\quad  Jinsong Liu\footnote{Shanghai University of Finance and Economics. Email: \texttt{liujinsong@163.sufe.edu.cn}}\quad \\ Jiyuan Tan\footnote{Fudan University.Email: \texttt{JiyuanTan19@gmail.com}} \quad Yinyu Ye\footnote{Standford University. Email: \texttt{yyye@stanford.edu}}}

\date{}
\begin{document}
\maketitle

\begin{abstract}
  SOLNP+ is a derivative-free solver for constrained nonlinear optimization. It starts from SOLNP proposed in 1989 by Ye \cite{yesolnp} with the main idea that uses finite difference to approximate the gradient. We incorporate the techniques of implicit filtering, new restart mechanism and modern quadratic programming solver into this new version with an ANSI C implementation. The algorithm exhibits a great advantage in running time and robustness under noise compared with the last version by MATLAB. SOLNP+ is free to download at \url{https://github.com/COPT-Public/SOLNP_plus}.
\end{abstract}
\section{Introduction}

	In many optimization problems, it's computationally difficult or even impossible to calculate the derivative of the objective or constraint functions, see \cite{audet2006finding}, \cite{booker1998managing}, \cite{booker1998optimization} for examples. Therefore, a class of nonlinear optimization methods not using derivative information, called derivative-free, BlackBox, gradient-free or zeroth-order optimization methods, have been proposed.

	We consider the following general constrained nonlinear optimization problem:
	\begin{equation}\label{problem}
	    \begin{aligned}
	        \min_{x} \quad &f(x) \\
            \text{s.t.}\quad &g(x) = 0\\
            &l_h \leq h(x) \leq u_h\\
            &l_x \leq x \leq u_x
	    \end{aligned}
	\end{equation}
    where $ f: \mathbb{R}^n\rightarrow \mathbb{R},g: \mathbb{R}^n\rightarrow \mathbb{R}^{m_1},h: \mathbb{R}^n\rightarrow \mathbb{R}^{m_2} $
	are smooth functions. 
	
	 In the new solver SOLNP+, we make significant improvements on both algorithm and implementation of SOLNP\cite{yesolnp}. SOLNP+ follows the framework of SOLNP and iteratively solves \eqref{problem} by linearizing the equality constraint and applying the augmented Lagrangian method in \cite{robinson1972quadratically} together with implicit filtering techniques. We use the finite difference to approximate the gradient. At the $k$-th iteration, we solve a linearly constrained optimization problem with an augmented Lagrangian objective function. We first check if $x^k$ is feasible for the linear equality constraints. If not, we solve a linear programming (LP) problem to 
	find an interior feasible solution. Then we use the BFGS technique to update the approximation of the Hessian matrix of the augmented Lagrangian objective function. After that we obtain sequential convex quadratic programmings (QP) to approximate the linearly constrained optimization problem with an augmented Lagrangian objective function. Then we solve sequential convex quadratic programmings starting from the feasible solution generated by LP.\par 
    Our main contributions can be summarized as follows: (1) SOLNP+ is able to deal with general nonlinear objective and constraint functions; (2) We exploit the adaptive step size of finite difference to approximate the gradient, which is more stable and relatively less sensitive to noise; (3) We provide a reliable, effective, open source, derivative-free solver. The numerical experiments show that it is one of the best solvers in its category.\par 
    Though some previous work considers the finite difference approach inefficient, our numerical results show that our algorithm is comparable with those methods that use trust-region with interpolation. Our observation is consistent with Berahas et al.\cite{berahas2019derivative}. Berahas et al. \cite{berahas2019derivative} argue that the high cost in implementing finite difference approach is often offset by its fast convergence speed.\par  
    
    \textbf{Outline of this Paper:} In section \ref{sec:review}, we review some literature in derivative-free optimization field and compare our algorithm with the algorithms in \cite{troltzsch2016sequential} and \cite{gilmore1995implicit}. In section \ref{sec:algo}, we reformulate problem (\ref{problem}) and introduce the ALM and SQP framework used by SOLNP+ in detail. Next, in section \ref{sec:details}, we describe some techniques we use in SOLNP+, including implicit filtering, restart and updating the penalty parameter. In section \ref{sec:numer}, we present the test results of SOLNP+ on three problems, the Hock and Schittkowski test problems \cite{hock1980test} with and without noise and the  pharmacodynamic problems \cite{arampatzis2019mu}.\par 
	\section{Literature Review}\label{sec:review}

	Over the past 30 years, there have been emerging demands in reality for optimizing problems whose derivatives are computationally expensive or even infeasible. Thus derivative-free methods come in handy. For example, in \cite{audet2006finding}, derivative-free optimization methods are used to tune parameters of nonlinear optimization methods. Direct search methods are used in \cite{meza1994direct} and \cite{alberto2004pattern} to solve molecular geometry optimization
	problems. In \cite{malik2019derivative}, Malik et al. study derivative-free methods for policy optimization over the class of linear policies.

	Numerous classic derivative-free optimization methods were first proposed to solve unconstrained problems, which are often adapted to solve constrained problems. There are mainly two classes of derivative-free unconstrained optimization methods in the literature.\par 
 
    \textbf{Direct Search Method:} Direct search methods take the next candidate point according to a finite number of sampled objective functions without derivative approximation or model building. Coordinate search(CS) \cite{fermi1952numerical}, 
	generalized pattern search(GPS) \cite{torczon1997convergence} and mesh adaptive directional search(MADS) \cite{audet2006mesh} are typical direct search methods as 
	they search along some directions to reduce the objective function value, while the differences among these algorithms lie in the set of alternative directions they choose. While Nelder Mead simplex method 
	\cite{nelder1965simplex} moves a simplex toward the optimal point by changing the vertices of the simplex instead of searching along some directions.\par 
 
	\textbf{Model-Based Method:} Trust region model-based methods exploit the trust region framework to build a local model that approximates the objective function in each iteration, which is relatively easier 
	to optimize. In \cite{powell1994direct}, the objective function is approximated by linear models. In \cite{powell2002uobyqa}, completely determined quadratic model is used to approximate the objective function, 
	which requires $(n+1)(n+2)/2$ function evaluations. While underdetermined quadratic models are exploited in \cite{powell2006newuoa} to reduce the number of function evaluations. 
	Using radial basis functions as models can approximate multimodal objective functions more precisely \cite{wild2008orbit}. There are also some hybrid methods. For example, 
	in \cite{kelley2011implicit}, Kelley uses the center difference to approximate gradient along with coordinate search, and BFGS models the hessian matrix to accelerate convergence.\par

	\textbf{Constrained Problems:} Derivative-free optimization methods deal with constraints in similar ways to those derivative-based methods. Thus the methods used in unconstrained cases can often be applied to solve constrained problems. 
	Generally, there are three approaches to dealing with constraints: using penalty function, filter and modeling constraints. As mentioned above, MADS is a directional direct search method, and along with extreme
	barrier \cite{audet2006mesh} or progressive barrier \cite{audet2009progressive}, it can also solve constrained cases. In \cite{troltzsch2016sequential}, nonlinear equality constraints are handled by 
	applying a sequential quadratic programming (SQP) approach. The Augmented Lagrangian method(ALM) is a general way to penalize constraints in nonlinear programming, where derivative-free methods can also exploit. 
	For example, In \cite{lewis2002globally}, Lewis and Torczon combined generalized pattern search and ALM framework. In \cite{audet2009progressive}, constraints are considered by a filter constructed from previously evaluated points.
	Noticing that in unconstrained cases, objective functions can be modeled as linear, quadratic or other functions, it is natural to model the constraints by a similar manner. In \cite{powell1994direct}, the constraints are modeled by linear interpolation, 
	which is penalized as a merit function to decide whether to accept the current point. Bajaj et al. \cite{bajaj2018trust} propose a two-phase algorithm to solve constrained problems. In the feasibility stage, the algorithm tries to find a feasible 
	point using modeled constraints. After finding such a point, the algorithm moves near the feasible region to minimize the objective function. \par 
	
     \textbf{Comparison with \cite{gilmore1995implicit} and \cite{troltzsch2016sequential}:} Both SOLNP+ and the implicit filtering method in \cite{gilmore1995implicit} use finite difference and BFGS update. However, the implicit filtering method can only deal with box constraints, while SOLNP+ can solve problems with general nonlinear constraints. \cite{troltzsch2016sequential} adopts the trust-region framework and uses quadratic polynomial interpolation to approximate the Lagrangian function and linear interpolation to approximate the constraints. It updates the dual variable using gradient ascend. SOLNP+ applies a modified augmented Lagrangian \cite{robinson1972quadratically} framework, which uses the dual solution of the QP subproblem as an approximation to the dual variable, and builds a quadratic model by BFGS update. \cite{troltzsch2016sequential} focuses on dealing with equality constraints, while SOLNP+ can handle both equality and inequality constraints.   
     
 
 	 \textbf{Improvements on SOLNP}: SOLNP was first proposed and implemented by Ye \cite{yesolnp}. SOLNP and SOLNP+ share the main iterative framework as mentioned above.
 	 The significant improvements of SOLNP+ can be summarized as follows:(1) When using the finite difference to approximate the gradient, we choose the step size adaptively to reduce the influence of local minimums caused by noise; (2)We perform the coordinate search at the end of the inner iteration to find a better initial point of the next inner iteration, which makes better use of calculated points; (3) We introduce the restart mechanism to avoid terminating at a suboptimal point and fix the inaccurate hessian matrices.\par 
\section{SOLNP+ Algorithm}\label{sec:algo}
In this section, we introduce the SOLNP+ algorithm in detail.

\subsection{Nonlinear Programming}
SOLNP+ solves the constrained nonlinear programming \eqref{problem}. 
Problem \eqref{problem} has strong modeling power since it includes various types of objective functions and constraints. From another perspective, however, its generality makes it an extremely hard problem.

By adding slacks to the inequality constraints, we convert the problem into
\begin{align*}
	\min\ & f(x) \\
	\text{s.t.}\ & g(x) = 0 \\
	& h(x) - s = 0 \\
	& l_{h} \leq s \leq u_{h} \\
	& l_{x} \leq x \leq u_{x}
\end{align*}
or for simplicity
\begin{equation} \label{p:sdnlp}
	\begin{aligned}
		\min\ & f(x) \\
		\text{s.t.}\ & g(x) = 0 \\
		& l_{x} \leq x \leq u_{x}
	\end{aligned}
\end{equation}
Here we abuse the notations for convenience. The following discussions are based on Problem \eqref{p:sdnlp}.

\subsection{Outer Iteration: Augmented Lagrangian Method}
SOLNP+ applies an Augmented Lagrangian Method (ALM) to solve \eqref{p:sdnlp}. At the $k$-th outer iteration of SOLNP+ with current solution $x^{k}$ and dual of the equality constraint $y^{k}$, we first approximate the equality constraints by its first-order Taylor expansion at $x^{k}$ and get
\begin{equation*}
	\begin{aligned}
		\min\ & f(x) \\
		\text{s.t.}\ & J^{k}(x - x^{k}) = -g(x^{k}) \\
		& l_{x} \leq x \leq u_{x}
	\end{aligned}
\end{equation*}
where $J^{k} = \left.\frac{\partial g}{\partial x}\right|_{x^{k}}$ is a numerical approximation to the Jacobian of $g(x)$. SOLNP+ does not require exact first-order information, so it is a zeroth-order method.

The modified augmented Lagrangian function in \cite{robinson1972quadratically} is
\begin{equation*}
	L_{\rho^{k}}^{k}(x,y) = f(x) - y^{T}\left[ g(x) - \left( g(x^{k}) + J^{k}(x - x^{k}) \right) \right] + \frac{\rho^{k}}{2}\left\| g(x) - \left( g(x^{k}) + J^{k}(x - x^{k}) \right) \right\|_{2}^{2}
\end{equation*}
Then we write the modified augmented Lagrangian problem
\begin{equation} \label{p:alm}
	\begin{aligned}
		\min\ & L_{\rho}^{k}(x,y^{k}) \\
		\text{s.t.}\ & J^{k}(x - x^{k}) = -g(x^{k}) \\
		& l_{x} \leq x \leq u_{x}
	\end{aligned}
\end{equation}
Notice that if we substitute the equality constraint of \eqref{p:alm} into $L_{\rho}^{k}(x,y^{k})$, the objective function becomes $f(x) - (y^{k})^{T}g(x) + \frac{\rho}{2}\|g(x)\|_{2}^{2}$, which is exactly the augmented Lagrangian function of \eqref{p:sdnlp} in classical ALM. However, we maintain the linearized equality constraints for better feasibility during iteration.

\subsection{Inner Iteration: Sequential Quadratic Approximation Approach}
Although we have linearized the hard constraints, \eqref{p:alm} is still a complicated problem since the objective function is a general nonlinear function. Therefore, we use BFGS to generate a series of symmetric positive definite matrices to approximate the Hessian of the objective function in \eqref{p:alm}. Then we obtain a sequence of convex quadratic programming (SQP) approximation problems with linear equality constraints.

Before solving the subproblem of SQP, we have to point out some important issues. We maintain the linearized equality constraints for better feasibility. However, the difference between the linearized constraints and those original ones sometimes unfortunately causes the SQP to be infeasible.

To fix the potential problem, we perform a two-step procedure: (1) We first find an interior (near-)feasible solution for SQP by projection or by solving linear programming (LP) and then adjust the constraint according to the feasible solution. (2) We perform an affine scaling algorithm (or some modern QP solvers as alternatives) to solve the SQP and regard the dual variable in SQP as the Lagrangian multiplier in \eqref{p:alm}.

\subsubsection{Find Feasible Starting Point}
We first need to find an interior feasible (or near-feasible) solution of \eqref{p:alm}. An Interior solution means that any component of the solution does not reach its boundary.
\paragraph{Free $x$ case}
If $x$ is free, we can get a feasible solution by projecting $x^{k}$ onto the hyperplane $J^{k}(x - x^{k}) = -g(x^{k})$, which is
\begin{equation*}
	x_{\text{feas}} = x^{k} - (J^{k})^{T}\left[ J^{k}(J^{k})^{T} \right]^{-1}g(x^{k}) 
\end{equation*}

\paragraph{Bounded $x$ case}
If $x$ is bounded, we try to solve the following interior linear programming \eqref{feas-LP} for an interior feasible (or near-feasible) solution.
\begin{equation}
	\begin{aligned}
		\min\ & \tau \\
		\text{s.t.}\ & J^{k}(x - x^{k}) - g(x^{k})\tau = -g(x^{k}) \\
		& l_{x} \leq x \leq u_{x} \\
		& 0 \leq \tau
	\end{aligned} \label{feas-LP}
\end{equation}
or in matrix form
\begin{align*}
	\min\ & c^{T}\begin{bmatrix} x \\ \tau \end{bmatrix} \\
	\text{s.t.}\ & A\begin{bmatrix} x \\ \tau \end{bmatrix} = b \\
	& l_{x} \leq x \leq u_{x} \\
	& 0 \leq \tau
\end{align*}
where $A = \begin{bmatrix} J^{k} & -g(x^{k}) \end{bmatrix}$, $b = J^{k}x^{k} - g(x^{k})$ and $c = (0,\cdots,0,1)^{T}$. Notice that $\begin{bmatrix} x^{k} \\ 1 \end{bmatrix}$ is an obvious feasible solution of \eqref{feas-LP} and if optimal $\tau = 0$, we find a feasible (but may be not interior feasible) solution of \eqref{p:alm}. Affine scaling method (see \cite{introLP} chapter 9 for more detail) can be applied to solve \eqref{feas-LP} to get the (near-)feasible solution $x^{k}_{0}$.

To ensure the following SQP is feasible, we need to change the right hand side of the equality constraint in \eqref{p:alm} to be $J^{k}(x^{k}_{0} - x^{k})$ instead of $-g(x^{k})$.

\subsubsection{Solve Sequential QP}
With $x^{k}_{0}$ and a guarantee of feasibility, we are now ready to solve \eqref{p:alm} by sequential quadratic approximation approach. SOLNP+ generates a series of convex QP problems by the BFGS method, which is
\begin{equation} \label{p:sqp}
	\begin{aligned}
		\min\ & \frac{1}{2}(x-x^{k}_{i})^{T} H^{k}_{i} (x-x^{k}_{i}) + \nabla L^{k}_{\rho}(x^{k}_{i},y^{k})^{T}(x-x^{k}_{i}) \\
		\text{s.t.}\ & J^{k} (x-x^{k}) = J^{k} (x^{k}_{0}-x^{k}) \\
		& l_{x} \leq x \leq u_{x}
	\end{aligned}
\end{equation}
and when $t^{T}s > 0$ we update $H^{k}_{i}$ as
\begin{equation} \label{u:bfgs}
	H^{k}_{i+1} = H^{k}_{i} + \frac{tt^{T}}{t^{T}s} - \frac{(H^{k}_{i}s)(H^{k}_{i}s)^{T}}{s^{T}H^{k}_{i}s}
\end{equation}
where $t = \nabla L^{k}_{\rho^{k}}(x^{k}_{i+1},y^{k}) - \nabla L^{k}_{\rho^{k}}(x^{k}_{i},y^{k})$ and $s = x^{k}_{i+1} - x^{k}_{i}$.

We can again apply affine scaling method \cite{ye1988interior} to solve \eqref{p:sqp} or use some modern QP solvers, e.g. OSQP \cite{osqp}. In the numerical experiments, sometimes we observe that the affine scaling method outperforms modern solvers in terms of the total function evaluation numbers. We consider that those ellipsoid constraints in affine scaling somehow work as trust regions and thus maintain the reliability of our quadratic approximation.

\subsubsection{Line Search}
To improve the solution $\tilde{x}_{i+1}^{k}$ obtained in \eqref{p:sqp} and meanwhile maintain feasibility, SOLNP+ uses bisection between $x^{k}_{i}$ and $\tilde{x}_{i+1}^{k}$ to get $x^{k}_{i+1}$ for lower objective value.

The line search approach has an influence on SOLNP+ from two sides. The positive aspect is that it may improve the solution quality and thus reduce the iteration number. The negative aspect is that the line search costs much function evaluation time, which is often expensive in zeroth-order applications. Therefore, we restrict the line search times in SOLNP+ to achieve a balance between the good and the bad.\par

\subsection{Summary of SOLNP+}
We have shown the whole SOLNP+ algorithm above. For a more clear presentation, we summarize SOLNP+ in Algorithm \ref{alg:SOLNP+}.

\begin{algorithm}[H]
	\caption{SOLNP+}
	\label{alg:SOLNP+}
	\begin{algorithmic}
		\INPUT $f,g,l_{x},u_{x}$ in problem, feasible starting point and its dual $(x^{0},y^{0})$, $\rho^{0}$ in ALM, inner and outer maximum iteration number $M, N$
		\OUTPUT Optimal solution and its dual $(x^{*}, y^{*})$
		\STATE /* --------------------- SOLNP+ in detail --------------------- */
		\FOR{$k = 0,\cdots, N$}
			\STATE /*  Sequential Quadratic Approximation Approach */
			\STATE Find an interior (near-)feasible solution $x^{k}_{0}$
			\STATE Adjust problem's feasibility using $x^{k}_{0}$
			\FOR{$i = 0,\cdots,M$}
				\STATE Solve QP \eqref{p:sqp} to get $\tilde{x}^{k}_{i+1}$
				\STATE Line search: Bisection between $x^{k}_{i}$ and $\tilde{x}^{k}_{i+1}$ to get $x^{k}_{i+1}$
				\STATE Update $H^{k}_{i}$ by BFGS as \eqref{u:bfgs}
				\IF{Inner stopping criteria satisfied}
					\STATE Break
				\ENDIF
			\ENDFOR
			\STATE 
			\STATE /* Update parameters and Convergence check */
			\STATE Update $\rho^{k}, H^{k}_{0}, y^{k}$
			\IF{Outer stopping criteria satisfied}
				\STATE Break
			\ENDIF
		\ENDFOR
	\end{algorithmic}
\end{algorithm}
\section{Implementation Details}\label{sec:details}
    In this section, we introduce some techniques we use to improve the solver's performance.

\subsection{Adaptively Choosing Gradient Step Size}\label{subsec:step}
    Choosing a proper step size when calculating the approximated gradient is significant when dealing with functions with noise. As the presence of noise may create a lot of local minimums, the algorithm will be stuck in these local minimums if it uses a small step size. To better motivate the strategy, we make the following calculation. Suppose that $ \delta $ is the step size and $ \epsilon $ is the noise, we calculate the gradient with finite difference,
    $$ [\nabla_\delta f(x)]_i = \underbrace{\dfrac{f(x+\delta e_i)-f(x)}{\delta}}_{(1) \text{Gradient estimation}} + \underbrace{\dfrac{\epsilon}{\delta}}_{(2)\text{Noise term}},$$
    where $ e_i = [0,\cdot,1,0,\cdots] $ is the standard basis of $ \mathbb{R}^n .$ If we take a small $ \delta $, the first term will have a good approximation to the gradient $ [\nabla f(x)]_i \approx (1) $. However, small $ \delta $ will make the second term larger, causing a larger error in the estimation. On the contrary, a large $ \delta $ can reduce the noise but cause some estimation error in the first term. There is a trade-off between the two terms.\par 
    We use the implicit filtering technique \cite{gilmore1995implicit} to deal with noise. The basic idea of \cite{gilmore1995implicit} is to adaptively choose $ \delta $ to implicitly filter out the noise of high frequency. We will change $ \delta^k_i $ at the end of each inner iteration dependent on how much the Lagrangian function has reduced. Let $ r^k_i = \dfrac{L^k_{\rho^k}(x^k_{i},y^k) - L^k_{\rho^k}(x^k_{i+1},y^k)}{\max\{1,L^k_{\rho^{k}}(x^k_{i},y^k)\}}. $ If $ r^k_i \geq c_{e}\delta^k_i$ , meaning that the Lagrangian function has reduced a lot, we will increase $ \delta_{i+1}^{k} = r_{ed} \delta_i^k $, expecting to further reduce it. If $ r^k_i \leq c_{\text{re}}\delta^k_i $, reduce $ \delta^{k+1}_1 = r_{rd} \delta_i^k .$ The inner iteration will be stopped if $\delta^k$ is modified.

\subsection{Combining Coordinate Search}
    The quadratic model constructed by the BFGS update may not be precise at the beginning of the optimization process. As a result, the model alone may not be able to reduce the function value. To make full use of the calculated points, we will record the best feasible point during calculation. 
    $$ \hat{x}_{i+1}^k = \text{argmin}_{|g(x+\delta^k_ie_i)|<\epsilon} f(x+\delta^k_ie_i). $$
    At the end of the inner iteration, we will compare $\hat{x}_{i+1}^k $ with the point obtained by the quadratic model and choose the better point as the initial point of the next inner iteration.  A similar approach has been considered in the implicit filtering algorithm \cite{gilmore1995implicit}.

\subsection{Restart Mechanism} \label{subsec:restart}
    SOLNP+ will restart in the following two situations.\par 
    First, if the relative difference between the two sequential outputs of the outer iteration is small, but our algorithm detects that the current point is far from the optimal point, we will perform a restart. In that situation, the algorithm may terminate at a suboptimal point if nothing is done. Formally, if 
    $$ \dfrac{f(x^{k-1})-f(x^{k})}{\max\{1,f(x^{k-1})\}} \leq \epsilon_s,$$ 
    where $ \epsilon_s $ is a small constant, the algorithm will check whether
    \begin{align}
        \left\Vert \mathcal{P}_{[l,u]}(x^*-\nabla_{\delta^k} L_k(x^k,y^k)) - x^k\right\Vert\leq \epsilon_a, \label{alm_stopapp}
    \end{align}
    where $ \epsilon_a $ is a large constant and $ \mathcal{P}_{[l_x,u_x]} $ is the projected operator. If (\ref{alm_stopapp}) is not satisfied, the algorithm will set the gradient step size to the initial value $ \delta = \delta_0 $ and modify Hessian $ H = \text{diag}(h^k_{11},\cdots,h^k_{nn}).$ Note that (\ref{alm_stopapp}) is an approximation to the standard optimal condition of the ALM method of the problem (\ref{p:sdnlp}),
    \begin{align}
        \mathcal{P}_{[l,u]}(x^*-\nabla L(x^*,y^*)) - x^* = 0.\label{alm_stop}
    \end{align}
    Using a larger $ \delta $ combined with coordinate search may take the solver out of some saddle points and further reduces the function value.\par 
    Second, if both the objective value and infeasibility increase after one outer iteration, which means the current estimated Hessian $H$  may be inaccurate  or step size  $\delta$  may be improper, the algorithm will set $y^k$ to zero vector and restart as the first case. 

\subsection{Update Penalty Parameter and Stop Criterion}
    We use the following heuristic to update the penalty parameter $\rho^k$. If the infeasibility at kth outer iteration $v^k$ satisfies $ v^k \leq c_z\cdot\text{tolerance}$, where$ c_z >1 $ is a constant, meaning that the infeasibility is small, we will set $ \rho = 0$. If the infeasibility increases a lot after one outer iteration, that is, $v^k \geq c_{ir} v^{k-1}$, we will increase $\rho^{k+1} = r_{ir} \rho^k, r_{ir}>1$. Otherwise, if  $v^k \leq c_{rr} v^{k-1}$, we will set  $\rho^{k+1} = r_{rr} \rho^k, 0<r_{rr}<1$.\par 
    The algorithm will stop the inner iteration if (1) the maximum inner iteration number has been reached or (2) the step size $\delta^k_i$ is improper, as mentioned in \ref{subsec:step}.  The outer iteration is ended if (\ref{alm_stopapp}) is satisfied, the tolerance of infeasibility is reached and the relative difference $\dfrac{f(x^{k-1})-f(x^{k})}{\max\{1,f(x^{k-1})\}}$ is small, meaning that the algorithm is unlikely to further reduce the objective value if continues optimizing. 
\section{Numerical Results}\label{sec:numer}
    In this section, we compare SOLNP+ with 2 famous derivative-free solvers: COBYLA \cite{powell1994direct} and NOMAD \cite{le2011algorithm}. We test them on classical benchmark sets \cite{hock1980test} and also compare them in a  pharmacodynamics problem.\par 


    COBYLA \cite{powell1994direct} is a derivative-free algorithm that uses linear polynomial interpolation to approximate the objective function and constraints. It is implemented in FORTRAN 77. In the experiment, we use the Matlab interface of PDFO \cite{ragonneau2021pdfo}, which is a cross-platform package providing interfaces for COBYLA.\par 

    NOMAD \cite{le2011algorithm} is a derivative-free solver that use the MADS algorithm \cite{audet2006mesh}. It is implemented in C++. In each step, it samples points to evaluate in the neighborhood of the current point to find a better one. It provides three ways to deal with general constraints: the extreme barrier, the progressive barrier and the filter. We use the progressive barrier in all the tests. As NOMAD can only deal with problems with inequality constraints, we have transformed the equality constraints in the tested examples into inequality in all tests. We use the Matlab interface of NOMAD version 3.9.1 in our experiment.\par 

    Apart from the initial point, SOLNP+ also needs an initial guess of the inequality constraints value. If the inequality constraint is bounded from two sides, i.e.
    $$ l \leq c(x) \leq u ,\quad \quad l,u \in \mathbb{R}^m,$$
    we will set the initial inequality constraints to $ i_0 = (l+u)/2.$ When the inequality constraints are bounded from one side, i.e., 
    $$ c(x) \leq u_x \quad \text{or} \quad c(x) \geq l_x, $$
    we use $i_0 = u_x - e  $ or $  l_x + e, e = [1,\cdots,1]^T $ as initial inequality constraints. We use affine scaling to solve the QP subproblems in all tests.\par 
    Although the three solvers use different stop criteria, the quality of solutions returned by the solvers is similar. All the experiments are tested on a laptop of Windows 10 with 2.60GHz 6-Core Intel i7-9750H processor and 16GB memory.\par 

\subsection{Functions without Noise}
    As SOLNP+ requires the initial point to lie in the interior of the box constraints, we choose all the problems in Hock and Schittkowski \cite{hock1980test} that satisfy this requirement to perform the test.\footnote{In experiments, we found that the initial and optimal values of problem hs54, hs70 and hs85 are inconsistent with the information provided by Hock and Schittkowski \cite{hock1980test}. Hence, we also exclude these problems.} It does not mean SOLNP+ cannot solve other problems, but changing the initial points may influence the solution of the problems. Hence, we only consider these 74 problems in our experiment.\par 
    The stopping tolerances of SOLNP+, NOMAD and COBYLA are set to $ 10^{-4} $. Hock and Schittkowski \cite{hock1980test} provide each problem with its solution. We call a problem solved if the solution $ x $  returned by the solver satisfies
    $$ f(x) - f_{\text{opt}} \leq 10^{-2}\max\{1,|f_{\text{opt}}|\} $$ 
    with infeasibility less than $ 10^{-4}$. \par
    
\begin{figure}[!htb]
        \centering
        \includegraphics{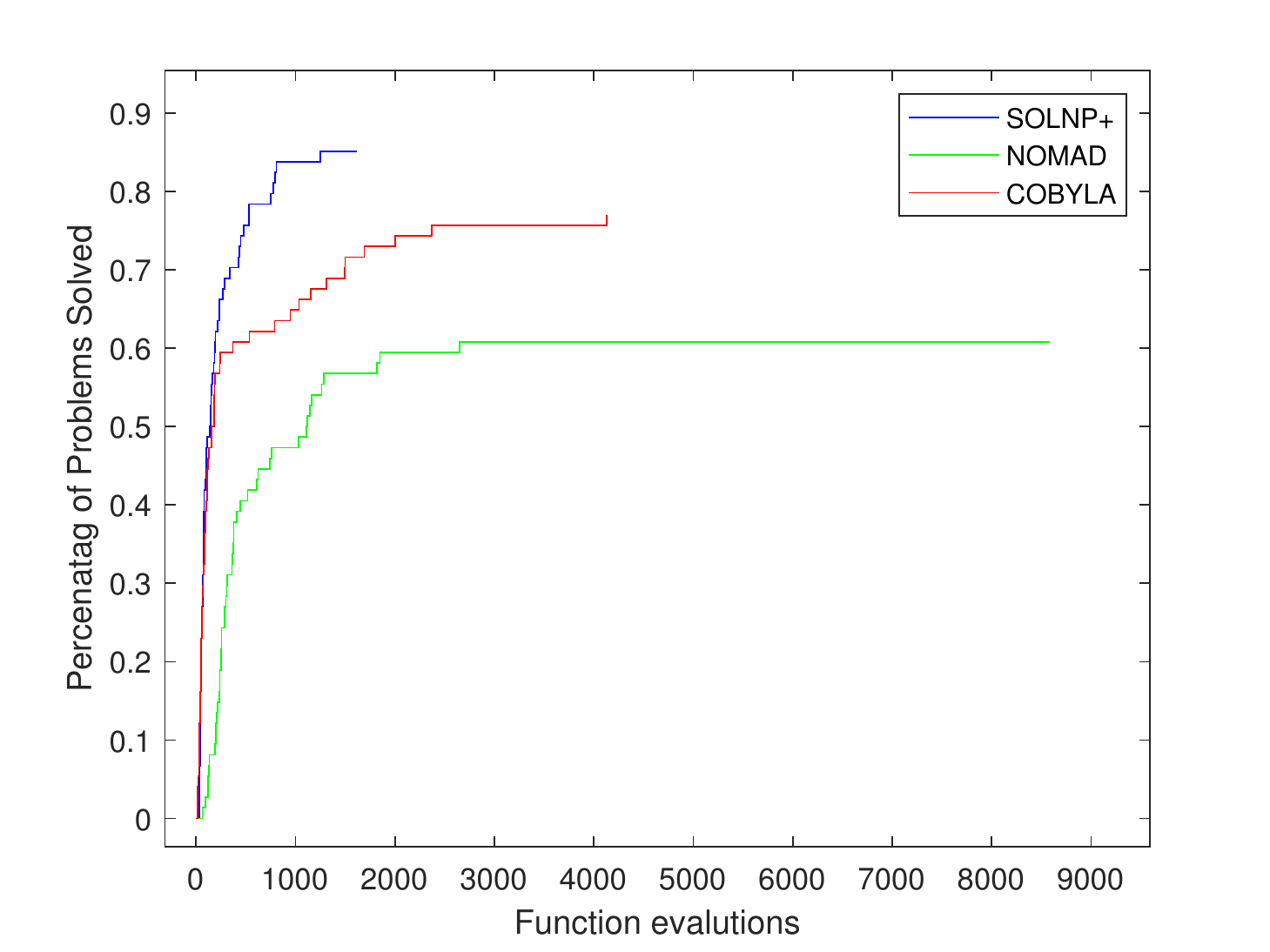}
        \caption{Test result of 74 problems in Hock and Schittkowski \cite{hock1980test} problems. Total running time of SOLNP+, NOMAD, COBYLA are 1.410250e+00s, 2.251209e+03s and 5.324220e+00s.}
        \label{fig:HS-fullprobs}
    \end{figure}

    Test results on 74 problems in Hock and Schittkowski \cite{hock1980test} are presented in Figure \ref{fig:HS-fullprobs}. We also present some problems with detailed information in Table \ref{exp:without_noise}. The reader can find the complete results in the appendix. The problems in Table \ref{exp:without_noise} are chosen arbitrarily from the 74 problems. They all have nonlinear constraints except hs38, which only has box constraints with a nonlinear objective function. \par 

    \begin{sidewaystable}[!htb]
        \small
        \center
        \begin{tabular}{|c|c|c|c|c|c|c|c|c|c|c|}\hline
        \multicolumn{1}{|c|}{\multirow{2}{*}{Prob.}}& \multicolumn{1}{|c|}{\multirow{2}{*}{Dim.}} & \multicolumn{3}{c|}{ Number of Evaluations} &\multicolumn{3}{c|}{Objective Function Value}&\multicolumn{3}{c|}{Running Time/s}\\ 
        \cline{3-11}& & SOLNP+ & NOMAD  & COBYLA & SOLNP+ & NOMAD &  COBYLA& SOLNP+ & NOMAD &  COBYLA \\ \hline 
        HS11 & 2 & 41  & 312 & 53 & -8.49787e+00  & \color{blue}-8.49846e+00 & \color{blue}-8.49846e+00
        &1.31810e-03 & 8.56036e+01 & 8.19570e-03\\ \cline{1-11}
        HS26 & 3 & 81 & 326 & 146 & {\color{blue}1.43427e-06}  & 3.56000e+00 & 2.11600e+01
        &1.83020e-03 & 2.65391e+01 & 1.41096e-02\\ \cline{1-11}
        HS38 & 4 & 165& 625 & 460  & 1.62759e-05 & \color{blue}2.25010e-13 & 7.87702e+00 
        &1.86670e-03 & 2.10570e+00 & 3.72709e-02 \\ \hline
        HS40 & 4 & 74 & 239 & 76  & \color{blue}-2.50025e-01 &-2.40655e-01 & -2.50000e-01
        &2.37800e-03 & 1.51294e+00 & 9.97500e-03 \\ \hline
        HS46 & 5 & 272 & 252 & 537  & \color{blue}4.30387e-09 &3.33763e+00 & 9.24220e-06
        &3.36330e-03 & 1.74984e+01 & 4.38302e-02 \\ \hline
        HS56 & 7 & 158  & 383  & 263  & -3.45603e+00 & -1.00000e+00& \color{blue}-3.45616e+00
        &4.07740e-03 & 1.46377e+01 & 2.52873e-02\\ \hline
        HS78 & 5 & 82 & 296 &110  &  \color{blue}-2.91974e+00 &2.73821e+00 &-2.91970e+00 
        &1.79030e-03 & 1.00057e+00 & 1.27184e-02\\ \hline
        HS79 & 5 & 75 & 353 &101  &  7.87804e-02 & 1.72669e-01& \color{blue}7.87768e-02
        &2.86990e-03 & 1.77057e+01 & 1.09047e-02\\ \hline
        HS80 & 5 & 104 & 312 & 96  &  \color{blue}5.39484e-02 &2.59025e-01 & 5.39498e-02
        &2.02920e-03 & 2.22013e+00 & 1.06562e-02 \\ \hline
        HS81 & 5 & 138& 328 &  153 &  \color{blue} 5.39470e-02&1.21224e-01 &  5.39498e-02
        &2.42640e-03 & 4.69137e+00 & 1.55220e-02 \\ \hline
        HS84 & 5 & 217 & 1818  &  54 &  \color{blue}-5.28034e+06 &-5.28019e+06 &-5.28033e+06 
        &7.55310e-03 & 1.36034e+01 & 7.56710e-03\\ \hline
        HS93 & 6 & 148& 1109 & 2367 &  1.35083e+02 &1.35525e+02 & \color{blue}1.35076e+02
        &5.06190e-03 & 1.01647e+01 & 2.09616e-01\\ \hline
        HS106& 8 & 530& 2670 & 4000 &  \color{blue}7.08435e+03 &7.66634e+03 & 8.94823e+03
        &1.21734e-02 & 4.79804e+01 & 3.87863e-01\\ \hline
        \end{tabular}
        \caption{Test results on Hock and Schittkowski \cite{hock1980test} problems. The blue color means that the solver returns an approximate optimal solution with better quality.}
        \label{exp:without_noise}
        \end{sidewaystable}
    \normalsize
    The numerical experiment shows that SOLNP+ solves more problems than the other two solvers on the test set. Besides, the SOLNP+ uses much fewer time for solving. From Table \ref{exp:without_noise} and appendix, we can see that SOLNP+ is much faster than the other two solvers in most problems, and the number of function evaluation is also comparable. 

\subsection{Functions with Noise}\label{sec:noise}
    
    We test our solver on functions with noises. To show that our algorithm is more robust under noise compared with the previous version, we also present the performance of SOLNP(the old version). 
    We add noise to the objective functions and constraints for the problem (\ref{problem}). That is, the solvers only observe $ \hat{f}(x) = (1+10^{-4}X(x))f(x), \hat{g}_i(x) = (1+10^{-4}Z_i(x))g(x), \hat{h}_i(x) = (1+10^{-4}N_i(x))h(x),$ where $ X(x),Z_i(x),N_i(x) \sim N(0,1), x\in \mathbb{R}^n $ are i.i.d. normal random variable. \par 
    The stopping tolerances of SOLNP+, NOMAD and COBYLA are set to $ 10^{-3} $. If the infeasibility of a solution is less than $ 10^{-3} $, we regard it as a feasible solution. The random seed in Matlab is set to $ 1 $ in all tests. The test results are presented in Table \ref{exp:with_noise}. For simplicity in presentation, we omit the detailed running time of each problem in the table.\par
    We can see from Table \ref{exp:with_noise} that except on HS26, HS38 and HS106, SOLNP+ keeps the good quality of solutions compared with the noiseless case, suggesting that it behaves well under some degree of noise.  On  HS78 and HS81, the relative difference between SOLNP+ and the best solver is less than $ 10^{-3} $. The previous version, SOLNP, is sensitive to the presence of noise. It fails to find a feasible point in most cases. Even if it finds a feasible point, the quality of the solutions worsens. The reason is that it uses a small step size when calculating gradients. The noise will cause large errors in gradient estimation. \par 
    COBYLA is influenced by the noise on HS46, HS56, HS81 and HS84. Compared with the noiseless case,  NOMAD is influenced by noise on HS79 and HS84. On HS84, NOMAD returns better feasible points with some failure. \par 
  
    \clearpage
\begin{sidewaystable}[!htb]
   \tiny
    \center
    \resizebox{1\textwidth}{0.25\textheight}{
    \begin{tabular}{|c|c|c|c|c|c|c|c|c|c|}\hline
    \multicolumn{1}{|c|}{\multirow{2}{*}{Prob.}}& \multicolumn{1}{|c|}{\multirow{2}{*}{Dim}} & \multicolumn{4}{c|}{ Average Number of Evaluations} &\multicolumn{4}{c|}{Average Objective Function Value}\\ 
    \cline{3-10}& & SOLNP & SOLNP+ & NOMAD  & COBYLA & SOLNP & SOLNP+ & NOMAD &  COBYLA \\ \hline 
    HS11 & 2  &118.13{\color{red}(20/50)} & 35.14 & 238.42& 43.54  &  4.03901e+03 &  -8.46861e+00 &\color{blue} -8.49988e+00& -8.42549e+00 \\ \hline 
    HS26 & 3  &125.55{\color{red}(21/50)} & 188.06{\color{red}(1/50)} & 213.24& 44.26  &  1.67703e+01 & \color{blue} 2.93551e-01 & 3.49606e+00& 2.11602e+01 \\ \hline 
    HS38 & 4  &37.00 & 224.44 & 702.12& 261.58  &  7.77777e+03 &  8.45308e-01 &\color{blue} 1.57504e-01& 7.93643e+00 \\ \hline 
    HS40 & 4  &512.17{\color{red}(44/50)} & 45.58 & 179.08& 67.14  &  -2.04409e-01 & \color{blue} -2.50324e-01 & -2.37238e-01& -2.49996e-01 \\ \hline 
    HS46 & 5  &127.00{\color{red}(28/50)} & 120.74 & 280.70& 101.02  &  2.76249e+00 & \color{blue} 4.44609e-05 & 3.33766e+00& 1.60209e+00 \\ \hline 
    HS56 & 7  &21.93{\color{red}(36/50)} & 531.78 & 377.60& 133.98 {\color{red}(1/50)} &  -1.00014e+00 &  -3.37944e+00 & -9.99998e-01&\color{blue} -3.45015e+00 \\ \hline 
    HS78 & 5  &--{\color{red}(50/50)} & 118.60 & 208.34& 73.58  &  -- &  -2.91860e+00 & -2.77044e+00&\color{blue} -2.91955e+00 \\ \hline 
    HS79 & 5  &889.00{\color{red}(47/50)} & 79.36 & 273.48& 79.62 {\color{red}(2/50)} &  3.75856e+00 &  7.88079e-02 & 4.27542e+01&\color{blue} 7.87840e-02 \\ \hline 
    HS80 & 5  &--{\color{red}(50/50)} & 87.18 & 221.14& 68.88  &  -- &  5.40269e-02 & 7.29409e-02&\color{blue} 5.39545e-02 \\ \hline 
    HS81 & 5  &1194.00{\color{red}(49/50)} & 141.74 & 223.58& 125.20 {\color{red}(1/50)} &  2.71448e-01 &  5.39633e-02 & 9.10489e-02&\color{blue} 5.39526e-02 \\ \hline 
    HS84 & 5  &17.96 & 236.44 & 589.86{\color{red}(36/50)}& 54.11 {\color{red}(41/50)} &  -2.35125e+06 &  -5.19516e+06 &\color{blue} -5.25703e+06& -5.24458e+06 \\ \hline 
    HS93 & 6  &19.00{\color{red}(39/50)} & 766.90 & 469.20& 86.38  &  1.37064e+02 &  1.36190e+02 &\color{blue} 1.35562e+02& 1.35922e+02 \\ \hline 
    HS106 & 8  &45.00{\color{red}(49/50)} & 581.98 & 1473.64& 82.30  &  1.49936e+04 &  1.50467e+04 &\color{blue} 7.80392e+03& 1.49971e+04 \\ \hline 
    
    \end{tabular}
    }
    \caption{Test results with noise on Hock and Schittkowski \cite{hock1980test} problems.  Each experiment is repeated 50 times. The blue color means that the solver returns a solution with better quality.``{\color{red}(fail time/total time)}'' means the number of times for which the solvers return an infeasible solution. The average is taken for all the feasible solutions returned by the solver. Total test time of SOLNP, SOLNP+, NOMAD and COBYLA are 3.49948e-01, 3.61491e-02, 1.25134e+02 and 1.34442e-01 seconds. }
    \label{exp:with_noise}
    \end{sidewaystable}
    \clearpage

\normalsize

\subsection{Application to Pharmacodynamics}
    We apply our algorithm to the tumor growth inhibition model \cite{ribba2012tumor}, which is solved with an evolution algorithm by Arampatzis et al. \cite{arampatzis2019mu}. The model aims to predict the tumor size and decide the best treatment plan for the patients. Ribba et al.\cite{ribba2012tumor} use the following ordinary differential equations to describe the growth of the tumor,
    \begin{align*}
        \dfrac{\d C}{\d t}  &= -\theta_1 C \\
        \dfrac{\d P}{\d t}  &= \theta_4P(1-\dfrac{P+Q+Q_P}{K}) +\theta_5Q_P-\theta_3P - \theta_1\theta_2CP \\
        \dfrac{\d Q}{\d t}  &= \theta_3P-\theta_1\theta_2CQ \\
        \dfrac{\d Q_P}{\d t}  &= \theta_1\theta_2CQ - \theta_5Q_P-\theta_6Q_P,
    \end{align*}
    with initial condition
    $$ C(0) = 0, \quad P(0) = \theta_7,\quad Q(0) = \theta_8, \quad Q_P(0) = 0, $$
    where $  \boldsymbol{\theta} = [\theta_1,\cdots,\theta_8] $ and $ K $ are constants. $ C(t) $ is drug concentration, and $ P, Q, Q_P $ are different tumor cells. The volume of tumor is defined to be $ P^* = P + Q + Q_P. $ The decision variable is $ t_1,\cdots,t_n,a_1,\cdots,a_n $. At the time $ t_i $, we give the patients drug at a dosage $ a_i $, and the differential equation is reinitiated by adding $ a_i $ to the drug concentration $ C(t_i) $ and keeping other variables the same. The goal of the model is to minimize the tumor size at the end of the treatment.\par
    
    However, the amount of drug need to satisfy some safety constraints \cite{harrold2009clinically}. First, the maximum amount of drug we administer one time should not exceed a scale,
    $$ 0 \leq a_i \leq 1,\quad \quad i = 1,\cdots n.$$
    Second, the maximum drug concentration cannot exceed the lethal dose,
    $$ 0 \leq \max_{t\in[0,t_{\text{end}}]} C(t) \leq v_{\max}.$$
    Third, the cumulative drug concentration should also be controlled,
    $$ 0\leq \int_0^{t_{\text{end}}}C(t) \d t \leq v_{\text{cum}}.$$ 
    Taking the three constraints into account, we can write the optimization problem as 
    \begin{align*}
        &\min_{t_1,\cdots,t_n,a_1,\cdots,a_n } P^* = P(t_{\text{end}})+Q(t_{\text{end}})+Q_P(t_{\text{end}})\\
        &\text{subject to} \quad\quad 0\leq t_i \leq t_{\text{end}}, \quad \quad i = 1,\cdots n,\\
        &\quad\quad\quad\quad\quad 0 \leq a_i \leq 1,\quad \quad i = 1,\cdots n,\\
        &\quad\quad\quad\quad\quad  0\leq \max_{t\in[0,t_{\text{end}}]} C(t) \leq v_{\max},\\
        &\quad\quad\quad\quad\quad 0\leq \int_0^{t_{\text{end}}}C(t)\d t \leq v_{\text{cum}}.
    \end{align*} 
    Notice that $ C(t) $ can be solved from the equation, i.e.,
     $$  C(t) = C(t_i)e^{-\theta_1(t-t_i)}, \quad \forall t \in [t_i,t_{i+1}] . $$
    We treat the calculation of $ P(t_{\text{end}}),Q(t_{\text{end}}),Q_P(t_{\text{end}}) $ as a blackbox to solve the problem.\par
\begin{figure}[H]
    
    \centering
    \includegraphics[width = 1\textwidth,height = 0.35\textheight,trim = 0 30 0 30]{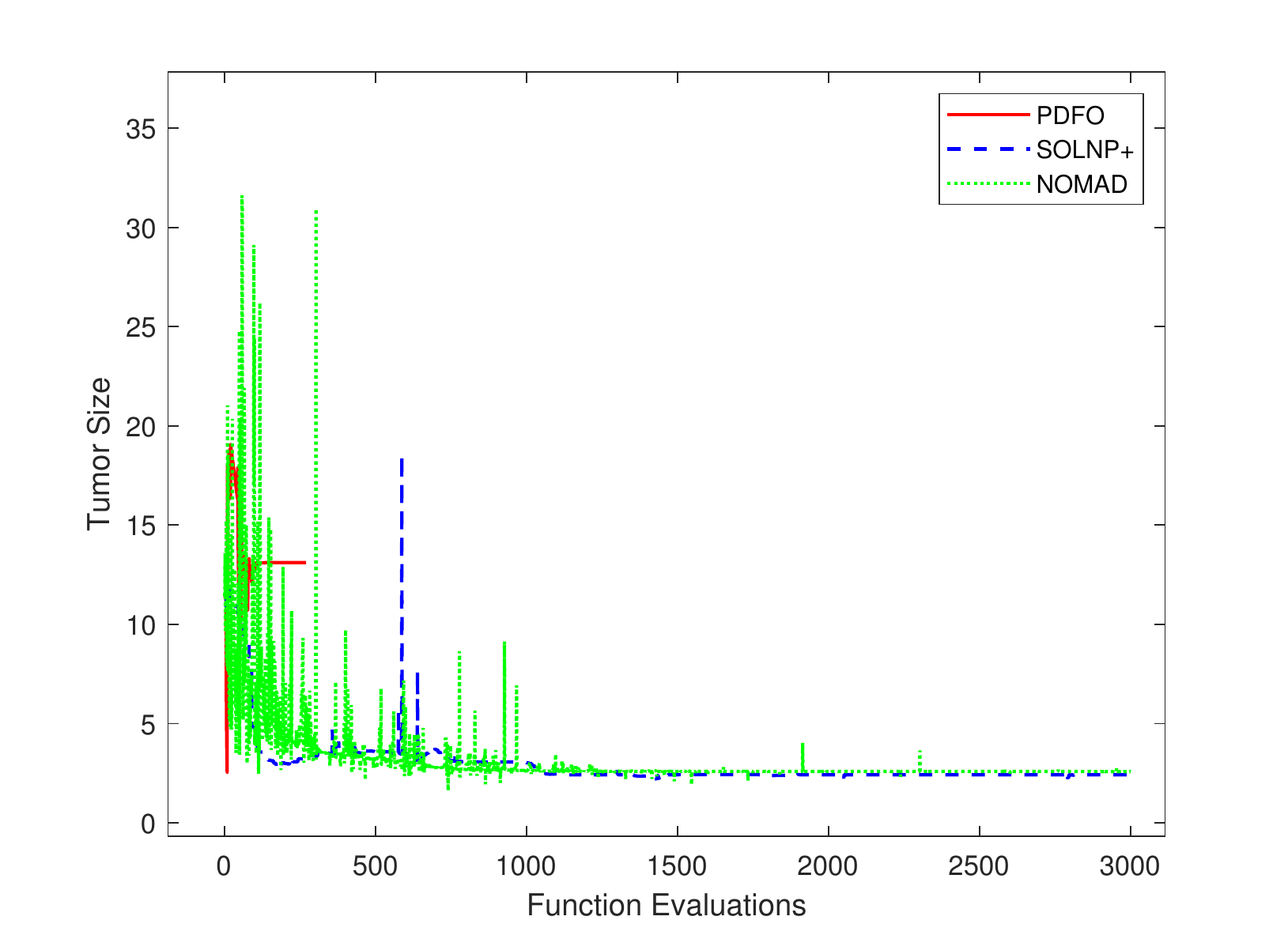}
    \caption{Convergence histories of the objective value. }
    \label{exp:tumor_his1}
\end{figure}

\begin{figure}[H]
    \centering
    \includegraphics[width = 1\textwidth,height = 0.35\textheight,trim = 0 30 0 30]{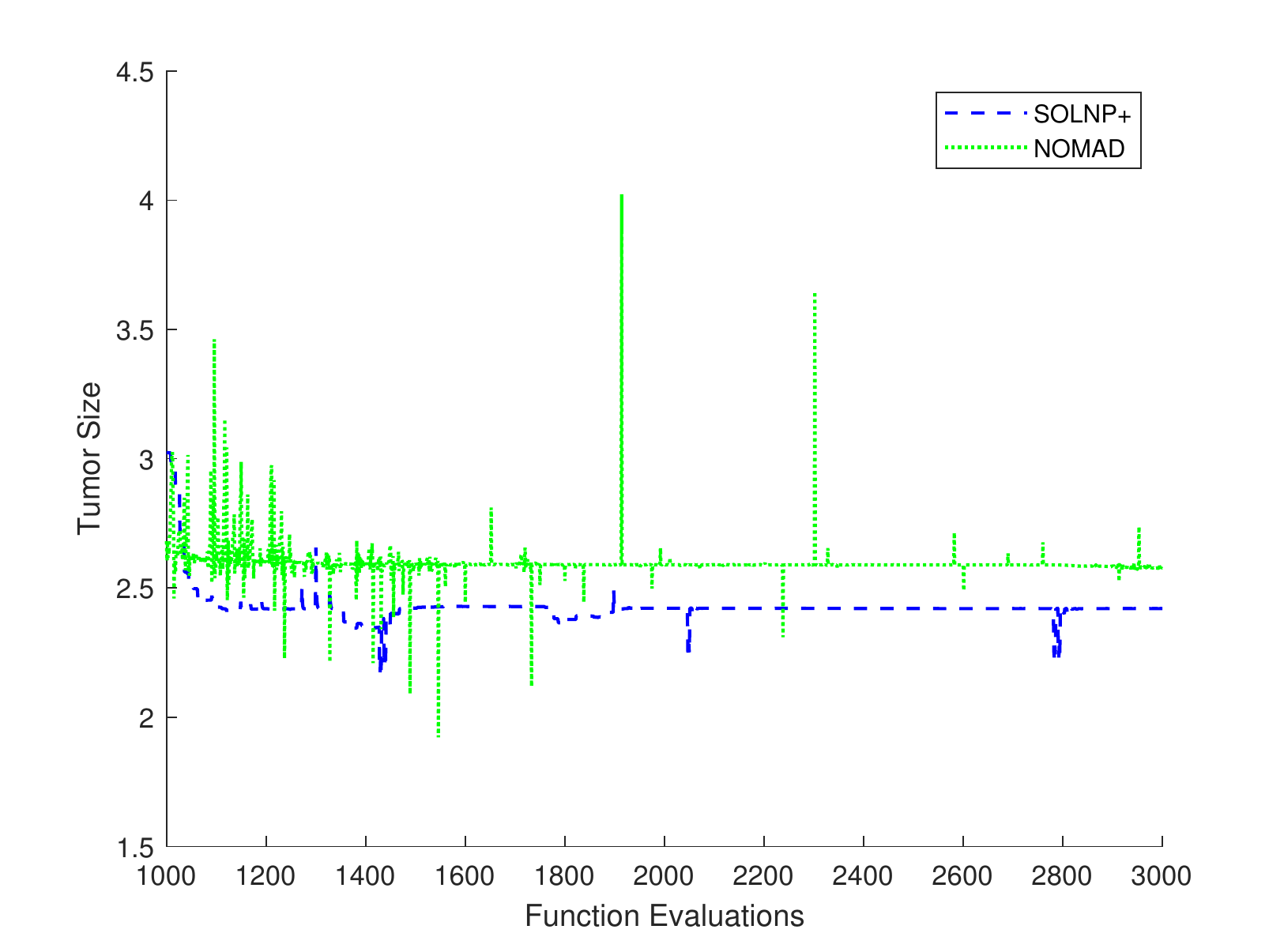}
    \caption{Convergence histories of the objective value after 1000 evaluations.}
    \label{exp:tumor_his2}
\end{figure}

\begin{figure}[H]
    \centering
    \includegraphics[width = 0.9\textwidth,height = 0.4\textheight,trim = 0 30 0 30]{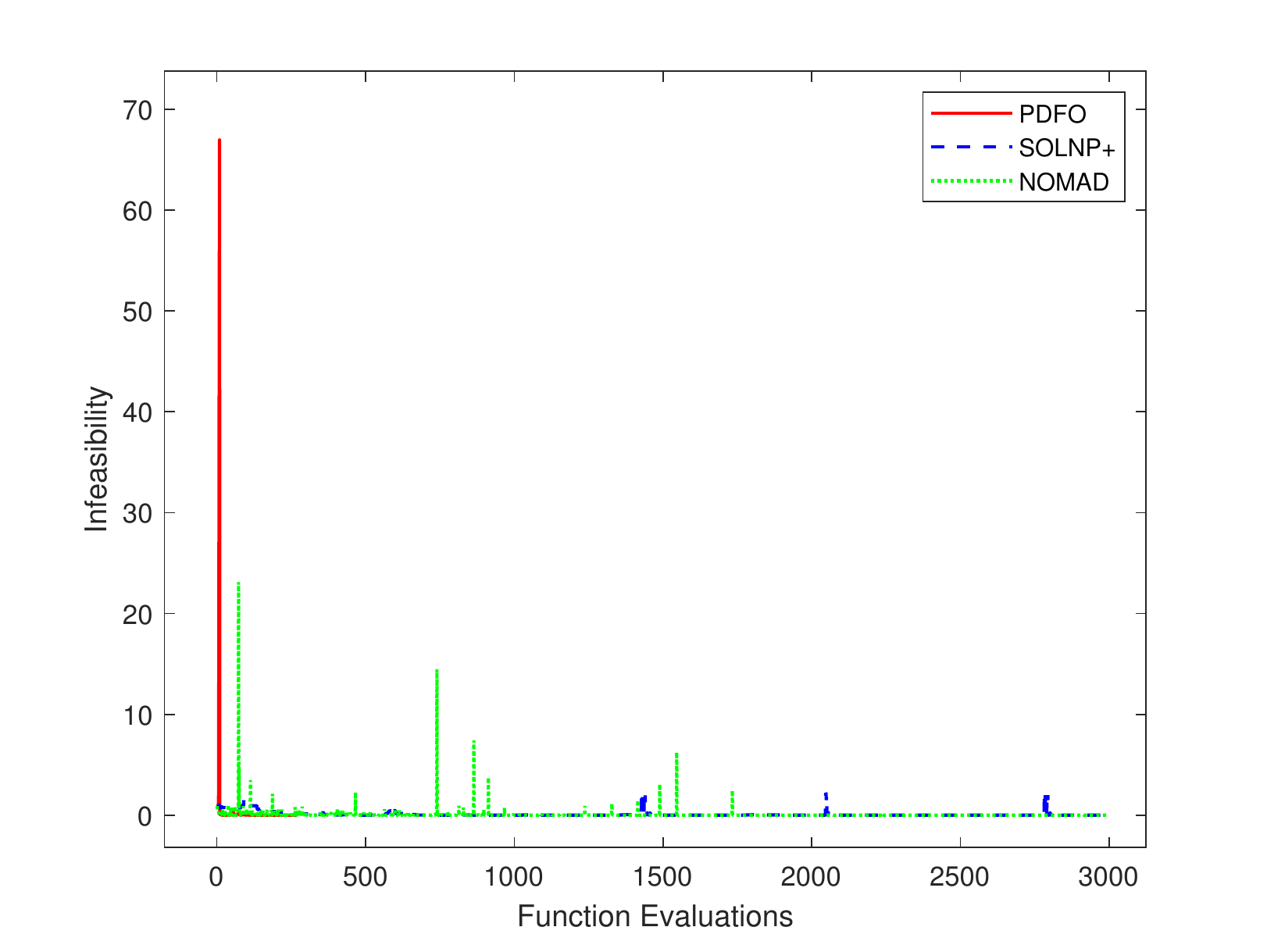}
    \caption{Convergence histories of the infeasibility.}
    \label{exp:tumor_his3}
\end{figure}

 
    Given the initial conditions, we calculate the approximate solution using the ode45 function in Matlab. In the experiment, we use the following constants, which are also used in \cite{arampatzis2019mu},
    $$ \boldsymbol{\theta} = [0.045, 4.52, 0.09, 0.11, 0.04, 0.00001, 0.09, 1],\quad K = 100. $$ 
    We consider $ n = 4,t_{\text{end}} = 200, v_{\max} = 1.1,v_{\text{cum}} = 65 $ with initial input as $ a_i = 0.5,t_i = t_{t_{\text{end}}}/2, i = 1,\cdots n$ (the default initial of SOLNP+, infeasible). The tolerance of solvers is set to $ 10^{-8}. $ and the maximum number of function evaluations is set to 3000. The test results are presented in Table \ref{exp:tumor_result} and Figure \ref{exp:tumor_his1}, \ref{exp:tumor_his2}, \ref{exp:tumor_his3}.\par

    \begin{table}[H]    
        \centering
        \begin{tabular}{|c|c|c|c|c|c|c|c|}\hline
        \multicolumn{1}{|c|}{\multirow{2}{*}{Problem}}& \multicolumn{1}{|c|}{\multirow{2}{*}{Dim}} & \multicolumn{3}{c|}{ Number of Evaluations} &\multicolumn{3}{c|}{Objective Function Value} \\ 
        \cline{3-8}& & SOLNP+ &NOMAD  & COBYLA & SOLNP+&NOMAD &  COBYLA  \\ \hline 
        Tumor & 8 & 3000 & 3000 &270 & 2.41949e+00& 2.57695e+00 & 1.31129e+01\\ \cline{1-8}
        \end{tabular}
        \\ [10pt]
        \centering
        \begin{tabular}{|c|c|c|c|c|c|}\hline
            \multicolumn{3}{|c|}{ Infeasibility} & \multicolumn{3}{c|}{Running Time/s}\\ 
            \hline  SOLNP+&NOMAD &  COBYLA & SOLNP+&NOMAD &  COBYLA \\ \hline 
            5.31037e-09& 0.00000e+00 &  0.00000e+00 &5.06375e+00 &3.28349e+01  &  9.16734e-01  \\ \hline
            \end{tabular}
        \caption{Final output in the tumor problem of three solvers.}
        \label{exp:tumor_result}
        \end{table}

We observe that SOLNP+ gets a solution with a better objective value among the three while keeping the infeasibility gap less than $ 10^{-8} $. We can see from the three figures that COBYLA stops after 270 function evaluations, ending with a feasible suboptimal solution. It finds a feasible point quickly in the beginning, but the objective value increases lately. On the contrary, SOLNP+ chooses to reduce the objective value even at the cost of some feasibility. After SOLNP+ finds a point with a small objective value, it chooses to reduce the infeasibility gap while keeping the objective value almost the same. The trajectory of NOMAD tends to oscillate dramatically at the beginning 1600 evaluations. We also notice that though both NOMAD and SOLNP+ make 3000 evaluations, the running time of SOLNP+ in this problem is significantly less.

\section{Summary}
    In this paper, we propose a new C implementation of the SOLNP+ solver for general constrained derivative-free optimization. We use the ALM framework and SQP approach to deal with the nonlinear constraints. The implicit filtering technique is used to increase the robustness of the algorithm under noise. The numerical experiments show that SOLNP+ is comparable with NOMAD and COBYLA. The solver is open source and still under active development. 
    

\bibliographystyle{plainnat}


\newpage
\appendix
\begin{sidewaystable}[!htb]
    \tiny
    \center

    \begin{tabular}{|c|c|c|c|c|c|c|c|c|c|c|c|}\hline
    \multicolumn{1}{|c|}{\multirow{2}{*}{Prob.}}& \multicolumn{1}{|c|}{\multirow{2}{*}{Dim.}}& \multicolumn{1}{|c|}{\multirow{2}{*}{Con.}}  & \multicolumn{3}{c|}{ Number of Evaluations} &\multicolumn{3}{c|}{Objective Function Value}&\multicolumn{3}{c|}{Running Time/s}\\\cline{4-12}
    & & & SOLNP+ & NOMAD  & COBYLA & SOLNP+ & NOMAD &  COBYLA& SOLNP+ & NOMAD &  COBYLA \\ \hline
    
    HS1 & 2 & 0  &112\color{red}(FAIL) & 134 & 1000\color{red}(FAIL) &1.46757e-01 & \color{blue}0.00000e+00 & 2.10833e+00&1.31900e-03 & 9.94425e-02 & 7.72052e-02 \\ \hline 
    HS3 & 2 & 0  &32 & 94 & 40 &1.01822e-03 & \color{blue}0.00000e+00 & 2.10605e-14&9.01400e-04 & 3.83821e-02 & 6.42740e-03 \\ \hline 
    HS4 & 2 & 0  &35 & 65 & 13 &2.66671e+00 & \color{blue}2.66667e+00 & \color{blue}2.66667e+00&8.29400e-04 & 3.02012e-02 & 6.01950e-03 \\ \hline 
    HS5 & 2 & 0  &44 & 120 & 33 &-1.91322e+00 & \color{blue}-1.91322e+00 & -1.91322e+00&1.09370e-03 & 8.05204e-02 & 6.05200e-03 \\ \hline 
    HS6 & 2 & 1  &67 & 175\color{red}(FAIL) & 40 &\color{blue}1.52385e-10 & 1.00000e+00 & 3.43024e-10&1.43740e-03 & 1.27828e+01 & 6.36010e-03 \\ \hline 
    HS7 & 2 & 1  &58 & 132 & 54 &\color{blue}-1.73207e+00 & -1.72918e+00 & -1.73205e+00&1.24500e-03 & 2.96294e+00 & 7.30550e-03 \\ \hline 
    HS8 & 2 & 2  &62 & 120 & 17\color{red}(FAIL) &\color{blue}-1.00000e+00 & \color{blue}-1.00000e+00 & -1.00000e+00&1.22270e-03 & 1.37847e-01 & 4.33970e-03 \\ \hline 
    HS9 & 2 & 1  &33 & 198 & 33 &-5.00000e-01 & \color{blue}-5.00000e-01 & -5.00000e-01&1.12690e-03 & 3.52492e+01 & 5.45970e-03 \\ \hline 
    HS10 & 2 & 1  &110 & 239 & 61 &-9.99946e-01 & \color{blue}-1.00000e+00 & -9.99998e-01&1.56150e-03 & 2.48453e+02 & 7.79450e-03 \\ \hline 
    HS11 & 2 & 1  &41 & 312 & 53 &-8.49787e+00 & \color{blue}-8.49846e+00 & -8.49846e+00&1.31810e-03 & 8.56036e+01 & 8.19570e-03 \\ \hline 
    HS12 & 2 & 1  &46 & 210 & 42 &-2.99942e+01 & \color{blue}-3.00000e+01 & -3.00000e+01&1.65880e-03 & 2.74775e+02 & 6.63450e-03 \\ \hline 
    HS14 & 2 & 2  &43 & 194 & 18\color{red}(FAIL) &\color{blue}1.39346e+00 & 1.39348e+00 & 2.29443e+00&1.94010e-03 & 1.67526e+01 & 5.25010e-03 \\ \hline 
    HS15 & 2 & 2  &540\color{red}(FAIL) & 218 & 124\color{red}(FAIL) &4.12014e+00 & \color{blue}3.06500e+02 & 3.60380e+02&7.30670e-03 & 2.88524e-01 & 1.69368e-02 \\ \hline 
    HS19 & 2 & 2  &52 & 378 & 28 &-6.96173e+03 & \color{blue}-6.96181e+03 & -6.96157e+03&1.97570e-03 & 1.08887e+01 & 5.74960e-03 \\ \hline 
    HS22 & 2 & 2  &37 & 256 & 17 &1.00004e+00 & \color{blue}1.00000e+00 & 1.00000e+00&1.77270e-03 & 3.63610e-01 & 4.73240e-03 \\ \hline 
    HS23 & 2 & 5  &196 & 232 & 18 &2.00013e+00 & \color{blue}2.00000e+00 & 2.00000e+00&5.56020e-03 & 4.85680e-01 & 4.82010e-03 \\ \hline 
    HS24 & 2 & 3  &232 & 255 & 14\color{red}(FAIL) &\color{blue}-1.00009e+00 & -9.99998e-01 & -2.13833e-02&5.24110e-03 & 1.33205e+00 & 4.75710e-03 \\ \hline 
    HS26 & 3 & 1  &81 & 326\color{red}(FAIL) & 146\color{red}(FAIL) &\color{blue}1.43427e-06 & 3.56000e+00 & 2.11600e+01&1.83020e-03 & 2.65391e+01 & 1.41096e-02 \\ \hline 
    HS27 & 3 & 1  &295\color{red}(FAIL) & 257 & 1500 &1.74806e-01 & \color{blue}4.00000e-02 & 4.09858e-02&2.76810e-03 & 1.21105e+01 & 1.15871e-01 \\ \hline 
    HS28 & 3 & 1  &61 & 363 & 60 &1.30568e-09 & \color{blue}0.00000e+00 & 2.98619e-09&1.36110e-03 & 2.15403e+01 & 7.99210e-03 \\ \hline 
    HS29 & 3 & 1  &72 & 445 & 61 &-2.26267e+01 & -2.26274e+01 & \color{blue}-2.26274e+01&1.30200e-03 & 1.00542e+02 & 9.34600e-03 \\ \hline 
    HS32 & 3 & 2  &41 & 332\color{red}(FAIL) & 21 &1.00003e+00 & 1.44418e+00 & \color{blue}1.00000e+00&1.54650e-03 & 1.24490e+01 & 4.53240e-03 \\ \hline 
    HS35 & 3 & 1  &42 & 250 & 51 &1.11133e-01 & \color{blue}1.11111e-01 & 1.11111e-01&2.02250e-03 & 3.43741e+01 & 7.50290e-03 \\ \hline 
    HS36 & 3 & 1  &41 & 362 & 31 &-3.29988e+03 & \color{blue}-3.30000e+03 & -3.29962e+03&1.47660e-03 & 5.91537e+00 & 5.86880e-03 \\ \hline 
    HS37 & 3 & 2  &48 & 290 & 72 &-3.45595e+03 & \color{blue}-3.45600e+03 & -3.45600e+03&1.92360e-03 & 1.64132e+02 & 9.21880e-03 \\ \hline 
    HS38 & 4 & 0  &165 & 625 & 460\color{red}(FAIL) &1.62759e-05 & \color{blue}2.25010e-13 & 7.87702e+00&1.86670e-03 & 2.10570e+00 & 3.72709e-02 \\ \hline 
    HS39 & 4 & 2  &146 & 300 & 115 &\color{blue}-1.00001e+00 & -1.00000e+00 & -1.00000e+00&2.54070e-03 & 6.74072e+00 & 1.26712e-02 \\ \hline 
    HS40 & 4 & 3  &74 & 239 & 76 &\color{blue}-2.50025e-01 & -2.40655e-01 & -2.50000e-01&2.37800e-03 & 1.51294e+00 & 9.97500e-03 \\ \hline 
    HS42 & 4 & 2  &93 & 309\color{red}(FAIL) & 94\color{red}(FAIL) &\color{blue}1.38577e+01 & 1.40000e+01 & 1.50000e+01&2.50750e-03 & 5.28478e+00 & 1.08098e-02 \\ \hline 
    HS43 & 4 & 3  &192 & 379 & 94 &-4.40000e+01 & \color{blue}-4.40000e+01 & -4.40000e+01&4.84940e-03 & 8.81464e-01 & 1.22515e-02 \\ \hline 
    HS46 & 5 & 2  &272 & 252\color{red}(FAIL) & 537 &\color{blue}4.30387e-09 & 3.33763e+00 & 9.24220e-06&3.36330e-03 & 1.74984e+01 & 4.38302e-02 \\ \hline 
    HS47 & 5 & 3  &1247 & 327\color{red}(FAIL) & 183 &\color{blue}-2.67155e-02 & 1.05254e+01 & -2.67142e-02&1.48170e-02 & 1.24320e+01 & 1.87143e-02 \\ \hline 
    HS48 & 5 & 2  &67 & 361\color{red}(FAIL) & 85 &\color{blue}6.48696e-15 & 6.76765e+01 & 1.73059e-09&2.07130e-03 & 1.79761e+01 & 8.87270e-03 \\ \hline 
    HS49 & 5 & 2  &103 & 539\color{red}(FAIL) & 1035 &\color{blue}5.34673e-09 & 2.70963e+01 & 3.22581e-05&1.95290e-03 & 1.95683e+01 & 8.12797e-02 \\ \hline 
    HS50 & 5 & 3  &276\color{red}(FAIL) & 367\color{red}(FAIL) & 182 &2.45549e+04 & 7.51600e+03 & \color{blue}2.28629e-08&3.07420e-03 & 7.04166e+00 & 1.62344e-02 \\ \hline 
    HS51 & 5 & 3  &73 & 761 & 88 &\color{blue}1.15959e-09 & 8.28525e-03 & 1.83869e-09&1.80280e-03 & 4.59296e+01 & 9.14510e-03 \\ \hline 
    HS52 & 5 & 3  &73 & 662\color{red}(FAIL) & 112 &\color{blue}5.32642e+00 & 6.77174e+02 & 5.32665e+00&1.53750e-03 & 1.16745e+01 & 1.18679e-02 \\ \hline 
    HS53 & 5 & 3  &42 & 512\color{red}(FAIL) & 89 &4.09304e+00 & 4.24501e+00 & \color{blue}4.09302e+00&1.46650e-03 & 5.91147e+00 & 9.91910e-03 \\ \hline 
    HS56 & 7 & 4  &158 & 383\color{red}(FAIL) & 263\color{red}(FAIL) &\color{blue}-3.45603e+00 & -1.00000e+00 & -3.45616e+00&4.07740e-03 & 1.46377e+01 & 2.52873e-02 \\ \hline 
    HS57 & 2 & 1  &34 & 285 & 32 &3.06464e-02 & \color{blue}2.84597e-02 & 3.06463e-02&1.17380e-03 & 4.08665e+01 & 5.35310e-03 \\ \hline

\end{tabular}
\caption{Test results on 74 Hock and Schittkowski \cite{hock1980test} problems. The blue color means that the solver returns an approximate optimal solution with better quality. {"\color{red}(FAIL)}" means either the infeasibility of the solution is larger than $10^{-4}$ or $f(x) - f_{\text{opt}} > 10^{-2}\max\{1,|f_{\text{opt}}|\}$, where $f_{\text{opt}}$ is the optimal solution.}
\end{sidewaystable}
\begin{sidewaystable}[!htb]
    \tiny
    \center

    \begin{tabular}{|c|c|c|c|c|c|c|c|c|c|c|c|}\hline
    \multicolumn{1}{|c|}{\multirow{2}{*}{Prob.}}& \multicolumn{1}{|c|}{\multirow{2}{*}{Dim.}}& \multicolumn{1}{|c|}{\multirow{2}{*}{Con.}}  & \multicolumn{3}{c|}{ Number of Evaluations} &\multicolumn{3}{c|}{Objective Function Value}&\multicolumn{3}{c|}{Running Time/s}\\\cline{4-12}
    & & & SOLNP+ & NOMAD  & COBYLA & SOLNP+ & NOMAD &  COBYLA& SOLNP+ & NOMAD &  COBYLA \\ \hline
    
    HS60 & 3 & 1  &55 & 238\color{red}(FAIL) & 53 &3.25691e-02 & 8.56727e+00 & \color{blue}3.25682e-02&1.72690e-03 & 1.11305e+01 & 6.96170e-03 \\ \hline 
    HS61 & 3 & 2  &428 & 230\color{red}(FAIL) & 73\color{red}(FAIL) &\color{blue}-1.43646e+02 & -9.90000e+01 & -8.19191e+01&3.70060e-03 & 1.31028e+01 & 7.96170e-03 \\ \hline 
    HS62 & 3 & 1  &80 & 374 & 180 &-2.62724e+04 & -2.62718e+04 & \color{blue}-2.62725e+04&1.91780e-03 & 1.03633e+01 & 1.68221e-02 \\ \hline 
    HS63 & 3 & 2  &186 & 198 & 66 &9.67518e+02 & 9.62123e+02 & \color{blue}9.61715e+02&2.39670e-03 & 2.63481e-01 & 8.11810e-03 \\ \hline 
    HS64 & 3 & 1  &382\color{red}(FAIL) & 1160 & 371 &6.37698e+03 & 6.29984e+03 & \color{blue}6.29984e+03&3.93140e-03 & 2.25190e+02 & 3.05943e-02 \\ \hline 
    HS67 & 3 & 14  &21\color{red}(FAIL) & 1029 & 1500\color{red}(FAIL) &-8.70720e+02 & \color{blue}-1.16212e+03 & -9.27492e+02&2.67510e-03 & 2.45412e+00 & 1.32240e-01 \\ \hline 
    HS68 & 4 & 2  &286 & 559\color{red}(FAIL) & 1313 &\color{blue}-9.20435e-01 & 2.40000e-05 & -9.20418e-01&1.32852e-02 & 8.60538e+00 & 2.02880e-01 \\ \hline 
    HS69 & 4 & 2  &178 & 340\color{red}(FAIL) & 2000 &\color{blue}-9.56713e+02 & -3.66467e+02 & -9.54907e+02&8.98280e-03 & 1.36046e+01 & 3.08403e-01 \\ \hline 
    HS72 & 4 & 2  &151 & 1262 & 947 &\color{blue}7.25129e+02 & 7.27739e+02 & 7.27679e+02&2.98710e-03 & 7.59217e+00 & 8.57387e-02 \\ \hline 
    HS73 & 4 & 3  &71 & 520 & 37\color{red}(FAIL) &\color{blue}2.98950e+01 & 3.01565e+01 & 3.04600e+01&3.47130e-03 & 6.14666e+00 & 6.01070e-03 \\ \hline 
    HS76 & 4 & 3  &56 & 1116 & 63 &-4.68178e+00 & -4.68181e+00 & \color{blue}-4.68182e+00&2.95350e-03 & 9.88028e+00 & 9.33300e-03 \\ \hline 
    HS77 & 5 & 2  &103 & 319\color{red}(FAIL) & 123 &\color{blue}2.41504e-01 & 1.17681e+01 & 2.41505e-01&2.11040e-03 & 9.16058e+00 & 1.22256e-02 \\ \hline 
    HS78 & 5 & 3  &82 & 296\color{red}(FAIL) & 110 &\color{blue}-2.91974e+00 & -2.73821e+00 & -2.91970e+00&1.79030e-03 & 1.00057e+00 & 1.27184e-02 \\ \hline 
    HS79 & 5 & 3  &75 & 353\color{red}(FAIL) & 101 &7.87804e-02 & 1.72669e-01 & \color{blue}7.87768e-02&2.86990e-03 & 1.77057e+01 & 1.09047e-02 \\ \hline 
    HS80 & 5 & 3  &104 & 312\color{red}(FAIL) & 96 &\color{blue}5.39484e-02 & 2.59025e-01 & 5.39498e-02&2.02920e-03 & 2.22013e+00 & 1.06562e-02 \\ \hline 
    HS81 & 5 & 3  &138 & 328\color{red}(FAIL) & 153 &\color{blue}5.39470e-02 & 1.21224e-01 & 5.39498e-02&2.42640e-03 & 4.69137e+00 & 1.55220e-02 \\ \hline 
    HS84 & 5 & 3  &217 & 1818 & 54 &\color{blue}-5.28034e+06 & -5.28019e+06 & -5.28033e+06&7.55310e-03 & 1.36034e+01 & 7.56710e-03 \\ \hline 
    HS88 & 2 & 1  &481 & 310 & 132 &\color{blue}1.28541e+00 & 1.36266e+00 & 1.36270e+00&1.45852e-01 & 7.60096e+00 & 1.01659e-01 \\ \hline 
    HS89 & 3 & 1  &338 & 408 & 57\color{red}(FAIL) &\color{blue}1.29181e+00 & 1.36266e+00 & 6.20616e-01&1.16437e-01 & 1.26777e+00 & 3.47683e-02 \\ \hline 
    HS90 & 4 & 1  &439 & 608 & 238 &\color{blue}1.28562e+00 & 1.36266e+00 & 1.36535e+00&1.58349e-01 & 1.79747e+01 & 2.10306e-01 \\ \hline 
    HS91 & 5 & 1  &777 & 744 & 1492 &\color{blue}1.29428e+00 & 1.36266e+00 & 1.36675e+00&2.63864e-01 & 3.41385e+01 & 9.99809e-01 \\ \hline 
    HS92 & 6 & 1  &797 & 1142 & 790 &\color{blue}1.29095e+00 & 1.36266e+00 & 1.36939e+00&3.85604e-01 & 1.68905e+01 & 8.18726e-01 \\ \hline 
    HS93 & 6 & 2  &148 & 1109 & 2367 &1.35082e+02 & 1.35525e+02 & \color{blue}1.35076e+02&5.06190e-03 & 1.01647e+01 & 2.09616e-01 \\ \hline 
    HS99 & 7 & 2  &532 & 432\color{red}(FAIL) & 870\color{red}(FAIL) &\color{blue}-8.31077e+08 & -8.11493e+08 & -8.31079e+08&7.74390e-03 & 2.03104e+01 & 8.00796e-02 \\ \hline 
    HS100 & 7 & 4  &448 & 1850 & 186 &6.80632e+02 & 6.80639e+02 & \color{blue}6.80631e+02&9.23220e-03 & 1.36042e+01 & 1.80305e-02 \\ \hline 
    HS104 & 8 & 5  &749 & 2651 & 1691 &\color{blue}3.95106e+00 & 3.95924e+00 & 3.95118e+00&1.77073e-02 & 1.85950e+00 & 1.71355e-01 \\ \hline 
    HS106 & 8 & 6  &530 & 2670\color{red}(FAIL) & 4000\color{red}(FAIL) &\color{blue}7.08435e+03 & 7.66634e+03 & 8.94823e+03&1.21734e-02 & 4.79804e+01 & 3.87863e-01 \\ \hline 
    HS108 & 9 & 13  &1136\color{red}(FAIL) & 1236\color{red}(FAIL) & 154 &-2.28470e+00 & -6.58523e-01 & \color{blue}-8.65929e-01&3.03705e-02 & 1.27335e+01 & 1.86281e-02 \\ \hline 
    HS110 & 10 & 0  &15\color{red}(FAIL) & 1282 & 193 &-4.31343e+01 & \color{blue}-4.57785e+01 & -4.57785e+01&1.08240e-03 & 6.58672e+00 & 1.82479e-02 \\ \hline 
    HS111 & 10 & 3  &236 & 776\color{red}(FAIL) & 4131 &-4.77441e+01 & -4.14032e+01 & \color{blue}-4.77611e+01&6.78570e-03 & 9.17814e+00 & 4.29204e-01 \\ \hline 
    HS112 & 10 & 3  &168\color{red}(FAIL) & 763\color{red}(FAIL) & 160\color{red}(FAIL) &-3.01912e+01 & -2.64018e+01 & -3.05820e+01&7.05950e-03 & 5.65460e+00 & 1.96422e-02 \\ \hline 
    HS113 & 10 & 8  &807 & 1620\color{red}(FAIL) & 242 &\color{blue}2.43066e+01 & 2.46032e+01 & 2.43069e+01&1.95854e-02 & 3.04938e+01 & 2.79286e-02 \\ \hline 
    HS114 & 10 & 11  &633\color{red}(FAIL) & 897\color{red}(FAIL) & 829\color{red}(FAIL) &-1.56320e+03 & -9.27265e+02 & -1.55294e+03&1.85496e-02 & 1.81414e+01 & 8.37127e-02 \\ \hline 
    HS117 & 15 & 5  &1618\color{red}(FAIL) & 8583\color{red}(FAIL) & 1153 &3.30784e+01 & 3.95038e+01 & \color{blue}3.23487e+01&3.67640e-02 & 3.61926e+02 & 1.38022e-01 \\ \hline

\end{tabular}
\caption{Test results on 74 Hock and Schittkowski \cite{hock1980test} problems. The blue color means that the solver returns an approximate optimal solution with better quality. {"\color{red}(FAIL)}" means either the infeasibility of the solution is larger than $10^{-4}$ or $f(x) - f_{\text{opt}} > 10^{-2}\max\{1,|f_{\text{opt}}|\}$, where $f_{\text{opt}}$ is the optimal solution.}
\end{sidewaystable}

\end{document}